\documentclass[11pt]{article}
\usepackage{amsmath,amsfonts,amsbsy,amsgen,amscd, amsxtra,amstext,amsopn,amssymb}
\usepackage{latexsym}
\usepackage{theorem}
\usepackage{mathrsfs}
\usepackage{pictex}
\usepackage{epic}
\usepackage{euscript}

\edef\qedrestoreat{\noexpand\catcode\lq\noexpand\@=\the\catcode\lq\@}%%
\catcode\lq\@=11

\ifx\protect\undefined\let\protect\relax\fi

\def\qed{\protect\@qed{$\qedsymbol$}}% suppressible printing of \qedsymbol
\def\pushright{\protect\@pushright}% takes an argument, 
\def\QED{\protect\@qed{{\rm Q.E.D.}}}% "quod   erat  demonstrandum"

\def\QEI{\protect\@qed{{\rm Q.E.I.}}}% "quod   erat  inveniendum"

\def\Proof{\protect\@Proof}\def\endProof{\protect\@endProof}%

\def\Proofof#1{\protect\@Proofof{#1}}\def\endProofof{\protect\@endProofof}%
\let\proofof\Proofof\let\endproofof\endProofof
\let\proofend\endproofof        
\def\qedsymbol{\raisebox{-.2ex}{$\Box$}}
\def\TheWordProof{\sc Proof.}
\def\TheWordProofof#1{\sc Proof of #1.}

\def\ProofFont{}

\newif\ifAutoQED\AutoQEDfalse

\newif\ifNumberResults
\ifx\theorem@style\undefined
   \NumberResultstrue
\fi

\def\parag@pushright#1{{% set up
    \parfillskip=0pt            % so \par doesnt push \square to left
    \widowpenalty=10000         % so we dont break the page before \square
    \displaywidowpenalty=10000  % ditto
    \finalhyphendemerits=0      % TeXbook exercise 14.32
    \hbox@pushright             % this used to be incorporated
    #1%                         % the end-of-proof mark (or whatever)
    \par}}%               % build paragraph with the above parameters

\def\hbox@pushright{% horizontal
    \unskip                     % remove previous space or glue
    \nobreak                    % don't break lines
    \hfil                       % ragged right if we spill over
    \penalty50                  % discouragement to do so
    \hskip.2em                  % ensure some space
    \null                       % anchor following \hfill
    \hfill                      % push \square to right
}%

\newif\if@qed\@qedfalse
\def\save@set@qed{\let\saved@ifqed\if@qed\global\@qedtrue}%
\def\restore@qed{\global\let\if@qed\saved@ifqed}

\def\@Proof{%
   \par\removelastskip\bigskip\penalty100
   \save@set@qed
   \noindent\ProofFont{\TheWordProof\enskip}%   
}%
\def\@Proofof#1{%
   \par\removelastskip\bigskip\penalty100
   \save@set@qed
   \noindent\ProofFont{\TheWordProofof{#1}\enskip}%
}%
\def\@endProof{%
   \qed\restore@qed
   \penalty-100 \medskip
}
\def\@endProofof{%
   \qed\restore@qed
   \penalty-100 \medskip
}
\def\@qed#1{%
\if@qed                                 % have we already done \qed?
     \global\@qedfalse%                - no - do it now, but not again
        \ifmmode\ifinner\pushright{#1}
        \else\eqno{\qedsymbol}\fi
        \else\pushright{#1}\fi%
\else\ifhmode\ifinner\else\par\fi\fi% -yes-just end paragraph, if any
\fi}
\def\@pushright#1{%
  {\ifvmode                       % vertical mode (see comments below)
       \null\hfill{#1}\par        % --- #1 alone in a paragraph!
  \else\ifmmode\maths@pushright{\hbox{#1}}% maths 
       \else\ifinner\hbox@pushright{#1}%  inside an \hbox
            \else\parag@pushright{#1}%    in a paragraph
  \fi  \fi  \fi
}}%
\def\maths@pushright#1{{%
  \ifinner
     \hbox@pushright{#1}%
  \else
     \eqno#1%   use TeX's right equation number feature (\eqno) 
     \def\]{$$\ignorespaces}% suppress LaTeX's error checking (HACK!)
  \fi
}}%
{\theorembodyfont{\slshape}
\newtheorem{theorem}{Theorem}[section]
\newtheorem{proposition}[theorem]{Proposition}
\newtheorem{lemma}[theorem]{Lemma}
\newtheorem{corollary}[theorem]{Corollary}
}
{\theorembodyfont{\rmfamily}
\newtheorem{definition}[theorem]{Definition}

\newtheorem{remark}[theorem]{Remark}
\newtheorem{example}[theorem]{Example}

}

\def\N{\mathbb{N}}

\def\R{\mathbb{R}}

\def\C{\mathbb{C}}

\def\P{\mathbb{P}}

\def\proj{\mathbb{P}}

\renewcommand{\a}{\alpha}

\renewcommand{\d}{\delta}

\newcommand{\e}{\varepsilon}
\newcommand{\g}{\gamma}

\newcommand{\s}{\sigma}

\renewcommand{\t}{\tau}

\newcommand{\z}{\zeta}

\def\scC{{\mathscr C}}
\def\scD{{\mathscr D}}

\def\mZ{{\mathcal Z}}

\newcommand{\x}{\times}
\newcommand{\<}{\langle}
\renewcommand{\>}{\rangle}

\renewcommand{\tilde}{\widetilde}

\def\grad{{\rm grad}\,}

\def\Oh{{\cal O}}

\def\algorithm{\begin{center}
               \begin{minipage}{6in}
               \begin{tabbing}
               \marks}
\def\falgorithm{\end{tabbing}
                \end{minipage}
                \end{center}}
\def\marks{nn\= nn\= nn\= nn\= nn\= nn\= nn\= \kill}

\def\II{{\mathrm{I\kern-0.5pt I}}}

\def\munorm{\mu_{\rm norm}}

\def\vol{{\mathsf{vol}}}
\def\rank{{\mathsf{rank}}}

\def\Tuproj{{T_{\P}^\perp}}
\def\Tusn{{T_{R}^\perp}}
\def\Tusplus{{T_{R}^+}}
\def\Tusminus{{T_{R}^-}}
\def\mc{{\mu}}
\def\kc{{\cal C}}
\def\nml{{\nu}}
\def\DD{{\mathcal D}}

\def\bfE{\mathop{\mathbf E}}
\def\Prob{\mathop{\rm Prob}}
\def\bfd{\mathbf d}
\def\Hd{{\mathcal H}_{\bfd}}

\def\cond{{\sf cond}}
\def\dist{{\sf dist}}

\def\munorm{\mu_{\rm norm}}
\def\munormR{\mu_{\rm norm,\R}\,}

\def\vol{{\mathsf{vol}}}
\def\rank{{\mathsf{rank}}}

\def\P{{\proj}}

\def\mc{{\mu}}
\def\kc{{\cal C}}
\def\tS{\Sigma_s}

\def\JoC{J. of Complexity}
\def\BAMS{Bull. Amer. Math. Soc.}
\def\JAMS{J. Amer. Math. Soc.}
\def\TCS{Theoretical Computer Science}

\begin{document}

\begin{title}
{\Large {\bf The probability that a small perturbation
of a numerical analysis problem is difficult}\thanks{Part of these results
were announced in C.R.\ Acad.\ Sci.\ Paris, Ser.~I 343 (2006) 145--150.}}
\end{title}
\author{Peter B\"urgisser\thanks{Dept.\ of Mathematics, University of
Paderborn, Germany. Partially supported
by DFG grant BU 1371.},
Felipe Cucker\thanks{Dept.\ of
Mathematics, City University of Hong Kong,
Kowloon Tong, Hong Kong. Partially supported by
CityU SRG grant 7001860.}, and Martin Lotz$^{\ddagger}$}
\date{}
\makeatletter
\maketitle
%\makeatother

\begin{quote}{\small
{\bf Abstract.}
We prove a general theorem providing smoothed analysis
estimates for conic condition numbers of problems of
numerical analysis.
Our probability estimates depend only on
geometric invariants of the corresponding sets of ill-posed inputs.
Several applications to linear and polynomial equation solving
show that the estimates obtained in this way
are easy to derive and quite accurate.
The main theorem is based on a volume estimate
of $\e$-tubular neighborhoods around a real algebraic subvariety
of a sphere, intersected with a disk of radius $\s$.
Besides $\e$ and $\s$, this bound depends only
the dimension of the sphere and on the degree of the defining equations.
}\end{quote}

%%%%%%%%%%%%%%%%%%%%%%%%%%%%%%%%%%%%%%%%

\section{Introduction}\label{se:intro}

In a seminal article~\cite{Demmel88} J.~Demmel suggested that
``to investigate the probability that a numerical analysis problem
is difficult, we need to do three things:
\begin{description}
\item{(1)}
Choose a measure of difficulty,
\item{(2)}
Choose a probability distribution on the set of problems,
\item{(3)}
Compute the distribution of the measure of difficulty induced
by the distribution on the set of problems.''
\end{description}
Then, for the measure of difficulty, Demmel proposed the {\em
condition number}. This is a positive number which, roughly speaking,
measures the sensitivity of the output to small perturbations of
the input. It depends only on the input data and the function being
computed. Condition numbers occur in endless instances of round-off
analysis.  They also appear as a parameter in complexity bounds for a
variety of iterative algorithms.

The main results in~\cite{Demmel88} carry out an analysis as sketched
in (1)--(3) above for the condition number $\scC$ of several problems.
This analysis exhibits bounds on the tail of the distribution of
$\scC(a)$,
showing that it is unlikely that $\scC(a)$ will be large.
>From these bounds one can obtain, using standard methods in probability
theory, bounds on the expected value of $\ln(\scC(a))$, estimating the
average loss of precision and average running time for algorithms
solving the considered problem. Demmel's results thus yield prime
instances of {\em average-case} analysis of algorithms in numerical
analysis.

While average-case analysis has undoubtedly advantages over worst-case
analysis, it is not itself without shortcomings, the most noticeable
being the arbitrariness of the selected probability distribution on the
set of inputs. To find a way out of these shortcomings,
D.~Spielman and S.-H.~Teng~\cite[\S3]{ST:02}
proposed a new form of analysis that arguably blends the best
of both worst-case and average-case. The idea is to replace showing that
\begin{quote}
``it is unlikely that $\scC(a)$ will be large''
\end{quote}
by showing that
\begin{quote}
``for all $a$ and all slight random perturbations $\Delta a$, it
is unlikely that $\scC(a+\Delta a)$ will be large.''
\end{quote}
A survey of this approach, called {\em smoothed analysis},
can be found in~\cite{ST:02,ST:06}.
If $\scD(c,\sigma)$ denotes a %symmetric
probability distribution centered at $c\in\R^{p+1}$ with
covariance matrix $\sigma^2{\rm id}_{p+1}$,
and {\bf E} denotes mathematical expectation,
we may summarize the objects of study of worst-case,
average-case, and smoothed analyses, for a function
$\psi:\R^{p+1}\to\R$, in the following table.
\bigskip

{\renewcommand{\arraystretch}{1.2}
\begin{tabular}[b]{||c|c|c||} \hline
{\tt worst-case analysis} &
{\tt average-case analysis} &
{\tt smoothed analysis} \\ [2pt] \hline
$\displaystyle \sup_{a\in\R^{p+1}} \psi(a)$ &
$\displaystyle \bfE_{a\in\scD(0,\sigma)} \psi(a)$ &
$\displaystyle \sup_{a\in\R^{p+1}}\bfE_{z\in\scD(a,\sigma)} \psi(z)$
\\ [10pt]\hline
\end{tabular}}
\bigskip

A remarkable feature of~\cite{Demmel88} is that the average-case
analysis performed there for a variety of problems is not done with
ad-hoc arguments adapted to the problem considered. Instead, these 
applications are all derived from a single result bounding the tail of the distribution
of $\scC(a)$ in terms of geometric invariants (degree and dimension)
of the set of ill-posed inputs of the problem for which $\scC$ is a
condition number. 

A first goal of this paper is to extend the results
of~\cite{Demmel88} from average-case to smoothed analysis. We will,
however, also prove average-case bounds. Demmel's paper dealt with
both complex and real problems. For complex problems he provided
complete proofs. For real problems, Demmel's bounds rely
on an unpublished (and apparently unavailable) result by A.~Ocneanu
on the volumes of tubes around real algebraic varieties. A second goal
of this paper is to prove a result akin to Ocneanu's (Theorem~\ref{th:realtubel}). 
We are not the first doing so. In~\cite{Wongkew:93}, R.~Wongkew gave 
a bound for the volume of tubes around real algebraic varieties. A number
of constants in his bounds, however, are not explicit and only
shown to be independent of the variety.

%%%%%%%%%%%%%%%%%%%%%%%%%%%%%%%%%%%%%%%%%%%%%%%%%%%%%%%%%%%%%%%
\subsection{Statement of the Main Result}

We assume our data space is $\R^{p+1}$, endowed with a
scalar product $\langle\ ,\ \rangle$.
By a semi-algebraic cone $\Sigma\subseteq \R^{p+1}$
we understand a semi-algebraic set $\Sigma\ne\{0\}$
that is closed by multiplications with positive scalars.
We say that $\scC$ is a {\em conic condition number} if there exists a
semi-algebraic cone $\Sigma\subseteq \R^{p+1}$, the set of {\em
ill-posed inputs}, such that, for all data $a\in\R^{p+1}\setminus\{0\}$,
\begin{equation*}
    \scC(a)=\frac{\|a\|}{\dist(a,\Sigma)},
\end{equation*}
where $\|\ \|$ and $\dist$ are the norm and distance
induced by $\langle\ ,\ \rangle$.

The best known condition number is that for
matrix inversion and linear equation solving. For a square matrix
$A$ it takes the form $\kappa(A)=\|A\|\|A^{-1}\|$
and was independently introduced by H.~Goldstine and
J.~von Neumann~\cite{vNGo47} and A.~Turing~\cite{Turing48}.
Strictly speaking, $\kappa(A)$ is not conic since the operator
norm $\|\ \|$ is not induced by a scalar product. Replacing
this norm by the Frobenius norm $\|\ \|_F$ yields the
(commonly considered) version
$\kappa_F(A) := \|A\|_F \|A^{-1}\|$ of $\kappa(A)$.
The Condition Number Theorem of C.~Eckart and G.~Young~\cite{EckYou}
then states that $\kappa_F(A)$ is conic, 
with $\Sigma$ the set of singular matrices.
Other examples
can be found in~\cite{ChC05}, where a certain property
(related with the so called level-2 condition numbers) is proved
for conic condition numbers. Furthermore, it is argued by 
Demmel in~\cite{Demmel87} that the condition numbers
for many problems can be bounded by conic ones.

Note that, since $\Sigma$ is a cone, for all $\lambda>0$,
$\scC(a)=\scC(\lambda a)$. Hence, we may restrict to
data $a$ lying in the sphere $S^p:=\{x\in\R^{p+1}\mid \|x\|=1\}$.
If we set $\tS:=\Sigma\cup(-\Sigma)$, then
the conic condition number $\scC$ can be estimated as
\begin{equation}\label{eq:conic-con}
  \scC(a)\leq \frac{\|a\|}{\dist(a,\tS)}
   = \frac{1}{d_{\P}(a,\tS\cap S^p)},
\end{equation}
where
%, abusing notation, $\Sigma$ is interpreted now as a subset of $S^p$ and
$d_{\P}$ denotes the projective distance in $S^p$,
which is defined as $d_{\P}(x,y)=\sin d_R(x,y)$ with $d_R$ being
the Riemannian (or angular) distance in $S^p$ (cf. Figure~1).
\begin{center}
   \input fig_real_1.pictex
\end{center}
Let $B_{\P}(a,\sigma)$ denote the open ball of radius $\arcsin\sigma$ 
around $a$ in $S^p$ ($\s$ corresponding to projective distance).  
We will endow $B_{\P}(a,\sigma)$ with the uniform probability measure.
Moreover, let
\begin{equation*}
 \Oh_{p} := \vol_p(S^p)
     = \frac{2\pi^{\frac{p+1}{2}}}{\Gamma\left(\frac{p+1}{2}\right)}
\end{equation*}
denote the $p$-dimensional volume of the sphere $S^p$.
Our main result is the following.

\begin{theorem}\label{th:main}
Let $\scC$ be a conic condition number with set of ill-posed inputs
$\Sigma$. Assume that $\Sigma\cap S^p\subseteq W$ where $W\subseteq S^{p}$
is the zero set in $S^p$ of homogeneous polynomials of degree at 
most~$d\ge 1$ and $W\ne S^{p}$.
Then, for all $a\in S^p$, all $\sigma\in(0,1]$, and all $t\geq 1$,
\begin{equation*}
   \Prob_{z\in B_{\P}(a,\sigma)}\{\scC(z)\geq t\}\leq
     4\sum_{k=1}^{p-1} {p\choose k} (2d)^k
  \left(1 +\frac{1}{t\s}\right)^{p-k} \left(\frac{1}{t\s}\right)^{k}
  +\, \frac{2p\Oh_p}{\Oh_{p-1}}\, (2d)^p \left(\frac{1}{t\s}\right)^{p}
\end{equation*}
and, for all $\sigma\in(0,1]$,
\begin{equation*}
   \sup_{a\in S^p}\,\bfE_{z\in B_{\P}(a,\sigma)}(\ln\scC(z))\leq
    2\ln p+2\ln d +2\ln\frac1\sigma+5.5.
\end{equation*}
In particular, for all $t\ge 1$ (take $\s=1$),
$$
   \Prob_{z\in S^p}\{\scC(z)\geq t\}\leq
         4\sum_{k=1}^{p-1} {p\choose k} (2d)^k\,
  \left(1 +\frac1{t}\right)^{p-k}\frac1{t^k}
  +\frac{2p\Oh_p}{\Oh_{p-1}}\,
  (2d)^p\,\frac1{t^p}
$$
and
$$
    \bfE_{z\in S^p}(\ln\scC(z))\leq
    2\ln p+2\ln d +5.5.
$$
\end{theorem}

The main idea towards the proof of Theorem~\ref{th:main}
is to reformulate
the probability distribution of a conic condition number
as a geometric problem in a sphere.
Indeed, for a measurable subset $V$ of $S^p$ 
we denote by $\vol_p(V)$ the $p$-dimensional volume of~$V$.
If $-V=V$ we define 
the {\em $\e$-neighborhood} around $V$ in $S^p$ by
$$
 T_{\P}(V,\e):= \{x\in S^p \mid d_{\P}(x,V) <\e\}.
$$
With this notation, using $\tS\cap S^p\subseteq W$,
we obtain from (\ref{eq:conic-con}) for $a\in S^p$ and  
$\sigma\in(0,1]$
\begin{eqnarray*}
  \Prob_{z\in B_{\P}(a,\sigma)}\left\{\scC(z)
  \geq\frac{1}{\e}\right\}
  &\leq&\Prob_{z\in B_{\P}(a,\sigma)}
      \left\{d_{\P}(z,\tS\cap S^p)\leq \e\right\}
  \leq\Prob_{z\in B_{\P}(a,\sigma)}\left\{d_{\P}(z,W)\leq \e\right\}\\
  &=&\frac{\vol_p(T_{\P}(W,\e)\cap B_{\P}(a,\sigma))}{\vol_p(B_{\P}(a,\sigma))}.
\end{eqnarray*}
The tail bounds in Theorem~\ref{th:main} will thus follow
from the following purely geometric statement.

\begin{theorem}\label{th:realtubel}
Let $W\subseteq S^p$ be a real algebraic variety
defined by homogeneous polynomials of degree at most $d\ge 1$
such that $W\ne S^p$.
Then we have for $a\in S^p$ and $0<\e,\s\le 1$
\begin{equation*}
 \frac{\vol_p \left(T_{\P}(W,\e)\cap B_{\P}(a,\s)\right)}
 {\vol_p B_{\P}(a,\s)} \le
  4\sum_{k=1}^{p-1} {p\choose k} (2d)^k\,
  \left(1 +\frac{\e}{\s}\right)^{p-k} \left(\frac{\e}{\s}\right)^{k}
  +\frac{2p\Oh_p}{\Oh_{p-1}}\,
  (2d)^p\,\left(\frac{\e}{\s}\right)^{p}.
\end{equation*}
In particular, (take $\s=1$)
$$
 \frac{\vol_p T_{\P}(W,\e)}{\Oh_p} \le
         4\sum_{k=1}^{p-1} {p\choose k} (2d)^k\,
  \left(1 +\e\right)^{p-k}\e^{k}
  +\frac{2p\Oh_p}{\Oh_{p-1}}\,
  (2d)^p\,\e^{p}.
$$
\end{theorem}

Here is a brief outline of the proof of Theorem~\ref{th:realtubel}:
The first step is an upper bound on the volume of an $\e$-neighborhood
of a smooth hypersurface in terms of integrals of absolute curvature
(Proposition~\ref{pro:tube-vol-mc}).
This is a variation of H.~Weyl's~\cite{weyl:39} exact formula
for the volume of tubes, a formula which, however, only holds
for sufficiently small~$\e$.
Then (Proposition~\ref{le:mc-estimate}) 
we derive a degree bound on these integrals of absolute curvature
based on the kinematic formula
of integral geometry and B\'ezout's theorem.
Finally we get rid of the smoothness assumption 
by some perturbation argument.

%%%%%%%%%%%%%%%%%%%%%%%%%%%%%%%%%%%%%%%%%%%%%%%%%%
%\subsection{Applications}\label{ssec:appl}

We will devote Section~\ref{sec:appl} to derive applications
of Theorem~\ref{th:main} to several condition numbers
occuring in the literature, namely, those for linear equation
solving, eigenvalue computation, polynomial system zero finding, 
and zero counting. 

\begin{remark}
Theorem~\ref{th:realtubel} could be stated for real projective
space $\P^p$ with the same bounds. While such a statement
is the most natural over the complex numbers
(cf.~\cite{BePa:05} and~\cite[Theorem~1.3]{BCL:06a}) it does not follow the 
tradition over the reals
(cf.~\cite{Demmel88,weyl:39,Wongkew:93}) and it is not a natural ambient
space for real conic condition numbers. Note that $\Sigma$ is not
necessarily symmetrical around the coordinate origin and that the
use of a symmetric $W$ (an algebraic cone containing the semi-algebraic
cone $\Sigma$) is just an artifice of our proofs.
\end{remark}

%%%%%%%%%%%%%%%%%%%%%%%%%%%%%%%%%%%%%%%%%%%%%%%%%%
\subsection{Relation to previous work}\label{sec:SA}

In most instances of smoothed analysis
(e.g.~\cite{CDW:05,DST,ST:02,ST:03,ST:04,Wsch:04})
one studies the behaviour, for a function $\psi:\R^{p+1}\to\R$,
of
\begin{equation}\label{eq:sa1}
   \sup_{a\in\R^{p+1}}\bfE_{z\in N^{p+1}(a,\sigma^2)} \psi(z)
\end{equation}
(possibly for sufficiently small $\sigma$) where $N^{p+1}(a,\sigma^2)$
denotes the $p+1$-dimensional Gaussian
distribution over $\R$ with mean $a$ and variance $\sigma^2$.
It has been argued that smoothed 
analysis interpolates between worst and average cases since it
amounts to the first for $\sigma=0$ and it approaches the second for
large $\sigma$.

When $\psi(\lambda x)=\psi(x)$ for all $\lambda>0$  ---e.g., 
a conic condition number--- it makes sense to restrict $\psi$ to 
the sphere $S^p$. In this case, it also makes sense to replace the
distribution $N^{p+1}(a,\sigma^2)$ by the uniform distribution
supported on the disk
$B_{\P}(a,\sigma)\subseteq S^p$
and to consider, instead of (\ref{eq:sa1}), the
following quantity
\begin{equation}\label{eq:sa2}
   \sup_{a\in S^p}\bfE_{z\in B_{\P}(a,\sigma)} \psi(z).
\end{equation}
Note that in this case, the interpolation mentioned above is
transparent. When $\sigma=0$ the expected value amounts to $\psi(a)$
and we obtain worst-case analysis, while if $\sigma=1$
the expected value is independent of $a$ and we obtain
average-case analysis.

It is this version of smoothed analysis, introduced
in~\cite{BCL:06a}, which we deal with in
this paper. Note that while, technically, this
``uniform smoothed analysis'' differs from the Gaussian
one considered so far, both share the viewpoint described
above.

We have already mentioned the
references~\cite{CDW:05,DST,ST:02,ST:03,ST:04,Wsch:04}
as instances of previous work in smoothed analysis.
In all these cases, an {\em ad hoc} argument is used to obtain
the desired bounds. This is in contrast with the goal
of this paper which is to provide general estimates
which can be applied to a large class of condition numbers.
We believe the applications in Section~\ref{sec:appl} give
substance to this goal.

We finish this section with a brief overview on previous work
on the relations between complexity, conditioning and probabilistic
analysis. In~\cite{Blum90}, L.~Blum suggested a complexity theory for
numerical algorithms parameterized by a condition number $\scC(a)$
for the input data (in addition to input size). 
S.~Smale~\cite[\S1]{Smale97} extended this suggestion by proposing to
obtain estimates on the probability distribution of~$\scC(a)$.
Combining both ideas, he argued, one can give probabilistic bounds
on the complexity of numerical algorithms.

The idea of reformulating probability distributions as quotients
of volumes in projective spaces (or spheres) to estimate
condition measures goes back at least
to Smale~\cite{Smale81} and Renegar~\cite{Ren87a}.
In particular, J.~Renegar~\cite{Ren87a} uses this idea to show bounds
on the probability distribution of a certain random variable in the
average-case analysis of the complexity of Newton's method.
Central to his argument is the
fact that this random variable can be bounded by a conic
condition number. The set of ill-posed inputs in~\cite{Ren87a}
is a hypersurface. An extension of these results to the case
of codimension greater than one was done by Demmel~\cite{Demmel88}
where, in addition, an average-case analysis of several
conic condition numbers is performed. Most of these results
are for problems over the complex numbers. An extension in another direction,
namely, to possibly singular ambient spaces, was done 
by C.~Beltr\'an and L.M.~Pardo~\cite{BePa:05}.
Another extension of Demmel's result, now to smoothed analysis for 
complex problems, was achieved in~\cite{BCL:06a}. 
\smallskip

The remaining of the paper is organized as follows. In 
Section~\ref{se:prelim} we provide the preliminary 
notations and results needed to prove Theorem~\ref{th:realtubel}. 
These are mostly taken from differential and integral geometry.
In Section~\ref{se:vol-tubes} we prove Theorem~\ref{th:realtubel} 
and Theorem~\ref{th:main}. 
Section~\ref{sec:appl} is devoted to several applications 
of our main result.
Finally, the precise value of some 
constants ---whose existence is well-documented in the literature but whose 
magnitude isn't--- is derived in an appendix.  

%%%%%%%%%%%%%%%%%%%%%%%%%%%%%%%%%%%%%%%
\section{Preliminaries}\label{se:prelim}
%\label{sec:machinery}

\subsection{Distances, volumes, and tubes on the sphere}
\label{se:dist-vol-tub}

The $p$-dimensional sphere $S^p$ carries the structure of a
compact Riemannian manifold.
Correspondingly, there is a Riemannian distance
$d_R(x,y)\in [0,\pi]$ between two points $x,y\in S^p$,
which is just the angle between these  points.
We denote by
$B_R(a,\a)=\{x\in S^p\mid d_R(a,x) < \a  \}$
the open ball of radius $\a$ centered at $a\in S^p$.

It will be more natural for us to work with the related notion of
{\em projective distance}, which is defined as
$d_\P(x,y) := \sin d_R(x,y) \in [0,1]$
(cf.~Figure~1).
We note that $d_\P$ satisfies all the axioms of a metric, except that
$d_\P(x,y)=0$ iff $x\in\{-y,y\}$.
In fact, $d_\P$ induces a metric on the real projective space $\P^p$
(obtained from $S^p$ by identifying antipodal points).
However, we prefer to work on the sphere, which seems more intuitive 
to us.

Let $V$ be a symmetric subset of $S^p$, that is $-V=V$.
For $0<\e\le 1$ we define the {\em $\e$-neighborhood} around $V$ by
$$
 T_{\P}(V,\e):= \{x\in S^p \mid d_{\P}(x,V) <\e\},
$$
where $d_{\P}(x,V) := \inf_{y\in V} d_\P(x,y)$.
This equals the $\a$-neighborhood of $V$ defined with respect to $d_R$,
where $\a=\arcsin\e$.
However note that
$T_\P(\{a\},\e)= B_R(a,\a) \cup B_R(-a,\a)$.
Somewhat inconsistently, we write
$B_\P(a,\e):=B_R(a,\a)$
and call this the open ball of radius $\e$ with respect to the
projective metric.

For a  measurable subset $A\subseteq S^p$ we write
$\vol_p A = \int_A dS^p$ for the {\em $p$-dimensional volume},
where $dS^p$ denotes the volume form induced by the Riemannian metric.
In order to compute volumes of balls and tubes around subspheres,
the following functions $J_{p,k}(\a)$ %defined for $1\le k\le p$ and $0\le\a\le\pi/2$
are relevant:
\begin{equation*}
  J_{p,k}(\a) := \int_0^\a (\sin\rho)^{k-1}\, (\cos\rho)^{p-k}\, d\rho \quad\mbox{($1\le k\le p$)}.
\end{equation*}

\begin{lemma}\label{le:J1}
For $1\le k\le p$, $0<\a\le\pi/2$, and $\e=\sin\a$, we have
$$
\vol_p T_{\P}(S^{p-k},\e) = \Oh_{p-k}\Oh_{k-1} J_{p,k}(\a),\quad
\vol_p B_R(a,\a)=\Oh_{p-1}J_{p,p}(\a).
$$
\end{lemma}

\begin{proof}
This follows from~\cite{weyl:39} or 
by straightforward calculation.
\end{proof}

The quantity $J_{p,k}(\a)$ can be easily bounded.
Recall that $\Oh_p$ denotes the $p$-dimensional volume of $S^p$. 

\begin{lemma}\label{le:J2}
The following estimates hold
(\/$1\le k\le p$, $0<\a\le\pi/2$, $\e=\sin\a$):
\begin{equation*}
 J_{p,k}(\a) \le \frac{\e^{k}}{k}\quad  \mbox{ if $k<p$},\quad
  \frac{\e^{p}}{p} \le J_{p,p}(\a) \le \frac{\Oh_p}{2\Oh_{p-1}}\, \e^p
\end{equation*}
with equality when $\a=\pi/2$ in the last upper bound.
\end{lemma}

\begin{proof}
To settle the first inequality note that for $k<p$
\begin{equation*}
  J_{p,k}(\a)
  \le \int_0^\a (\sin\rho)^{k-1}\, (\cos\rho)\, d\rho
   = \int_0^{\e} u^{k-1}\, du  = \frac{\e^{k}}{k}.
\end{equation*}
Similarly,
\begin{equation*}
 J_{p,p}(\a) = \int_0^{\a} (\sin\rho)^{p-1}\, d\rho \ge
   \int_0^{\a} (\sin\rho)^{p-1} \cos\rho\, d\rho = \frac{\e^p}{p}.
\end{equation*}
It is easy to check that $\a\to J_{p,p}(\a) (\sin\a)^{-p}$ is
monotonically increasing on $[0,\pi/2]$ by computing the derivative
of this function. Hence,
$J_p(\a) (\sin\a)^{-p} \le J_{p,p}(\pi/2)$.
% = \frac{\Oh_p}{2\Oh_{p-1}}.
>From Lemma~\ref{le:J1} we get
$ \frac12\Oh_p=\vol_p B_R(a,\pi/2)=\Oh_{p-1}J_{p,p}(\pi/2)$
from which it follows that $J_{p,p}(\pi/2)=\frac{\Oh_p}{2\Oh_{p-1}}$.
\end{proof}

In this paper, the notions of manifold and differentiability always 
refer to $C^\infty$-differentiability.
For a submanifold $M$ of $S^p$ and $0<\a\le\pi/2$ we define the
{\em $\a$-tube} $\Tusn(M,\a)$ around~$M$ by (compare with~\cite[p.~34]{gray:90})
\begin{align*}
 \Tusn(M,\a) := \{x\in S^p \mid\;& \mbox{there is a great circle segment in
                 $S^p$ of length $<\a$}\\
  & \mbox{from $x$ to $M$ that intersects $M$ orthogonally}\}.
\end{align*}
%\begin{align*}
% \Tuproj(M,\e) := \{x\in S^p \mid\;& \mbox{there is a line segment in
%      $S^p$ of projective length $<\e$}\\
%  & \mbox{from $x$ to $M$ that intersects $M$ orthogonally}\}.
%\end{align*}
Here we used the Riemannian distance. 
Sometimes, when thinking in terms of the projective distance 
and $M=-M$, 
it will be convient to use the notation of 
{\em $\e$-tube} $\Tuproj(M,\e):=\Tusn(M,\arcsin\e)$ defined  
for $0<\e\le 1$. 
Clearly, $\Tuproj(M,\e)$ can be characterized in a way similar to 
$\Tusn(M,\a)$. 
We note that $\Tuproj(M,\e)\subseteq T_{\P}(M,\e)$ and the inclusion is
in general strict. It can be shown, however, that if $M$ is a compact
submanifold, then equality holds.

We also need the notion of the $m$-dimensional volume
(Hausdorff measure) of subsets $T$ of $S^p$. For simplicity,
we restrict ourselves to semialgebraic sets, cf.~\cite{bocr:98}.
Let $T$ be an $m$-dimensional semialgebraic subset of $S^p$.
The Zariski closure $W$ of~$T$ in $S^p$ is a real algebraic variety
of dimension~$m$ and its regular locus $\mathrm{Reg}(W)$ is an 
$m$-dimensional submanifold of $S^p$.
We define the $m$-dimensional volume of $T$ by
$\vol_m T := \vol_m (T\cap \mathrm{Reg}(W))$.
This makes sense since $T\setminus\mathrm{Reg}(W)$ has dimension strictly
less than~$m$. For $k<m$ we set $\vol_k T:=0$.

%%%%%%%%%%%%%%%%%%%%%%%%%%%%%%%%%%%%
\subsection{A useful transformation formula}\label{se:trafo}

We will repeatedly use the following special case of the coarea formula.
A proof can be found in~\cite[III\S2, Folgerung~1]{suwi:72} 
or~\cite[Appendix]{howa:93}.

\begin{proposition}\label{pro:trafo}
Let $M,N$ be Riemannian manifolds of the same dimension and
$\varphi\colon M\to N$ be differentiable.
Suppose that $\int_M |\det D\varphi|\, dM$ is finite.
Then the fiber $\varphi^{-1}(y)$ is finite for almost all $y\in N$
and we have
$$
 \int_M |\det D\varphi|\, dM = \int_{y\in N} |\varphi^{-1}(y)|\, dN(y) .
$$
\end{proposition}

%%%%%%%%%%%%%%%%%%%%%%%%%%%%%%%%%%%%%%%%%%%%%%%%%
\subsection{Some differential geometry of hypersurfaces on spheres}
\label{se:dg-sphere}

For the following material from differential geometry we refer to
\cite{spiv3:79,thor:94}.

In the following let $M$ be a compact oriented smooth hypersurface of $S^p$
interpreted as a Riemannian submanifold. The orientation corresponds
to the choice of a unit normal vector field $\nml\colon M\to \R^{p+1}$ on $M$. 
The {\em Weingarten map} $L_M(x)\colon T_xM\rightarrow T_xM$ of $M$ at $x$ is
the linear map defined by $L_M(x):=-D\nml(x)$ (it is easy to verify
that this is a well-defined map). The {\em second
fundamental form} of $M$ at $x\in M$ is the corresponding bilinear
map $\II_M(x)\colon T_xM\times T_xM\to\R$, defined by $\II_M(x)(Y,Z)=
\langle L_M(x)Y,Z\rangle$ for all $Y,Z\in T_xM$. We are going to
describe these notions in terms of local coordinates of~$M$. Thus
let $v=(v_1,\ldots,v_{p-1})\mapsto x(v_1,\ldots,v_{p-1})\in
M\subset\R^{p+1}$ be a local parametrization of $M$. Then
\begin{equation*}
 \II_M(x)\left(\frac{\partial}{\partial v_i},\frac{\partial}{\partial v_j}\right)
    =-\left\langle \frac{\partial \nml}{\partial v_i},
     \frac{\partial x}{\partial v_j}\right\rangle
    = \left\langle\nml,\frac{\partial^2 x}{\partial v_i\partial v_j}\right\rangle,
\end{equation*}
where the last equality follows from deriving $\<\partial x/\partial
v_j,\nml\>=0$. In particular, $\II_M(x)$ and $L_M(x)$ are symmetric.
The eigenvalues $\kappa_1(x),\ldots,\kappa_{p-1}(x)$ of $L_M(x)$ are
called {\em principal curvatures} at~$x$ of the hypersurface~$M$.

\begin{example}\label{ex:sphere}
Consider the case of $M=S^{p-1}$, the subsphere of $S^p$ given by
the equation $x_p=0$. Then it easy to see that $\II_M(x)=0$ for all
$x\in M$. Hence all the principal curvatures of $M$ are zero. This
example makes clear that the principal curvatures are relative to
the ambient space~$S^p$. (Of course, $S^{p-1}$ is curved; however,
its ``curvature relative to the ambient sphere'' is zero.)
\end{example}

For $1\le i < p$ we define the {\em $i$th curvature} $K_{M,i}(x)$
of $M$ at $x$ as the $i$th elementary symmetric polynomial in
$\kappa_1(x),\ldots,\kappa_{p-1}(x)$, and
put $K_{M,0}(x):=1$.
In particular, $K_{M,p-1}(x) = \det L_M(x)$. Note that the $i$th
curvatures are essentially the coefficients of the characteristic
polynomial of the Weingarten map:
\begin{equation}\label{eq:charpol}
 \det(\mathrm{id}_{p-1} + t L_M(x)) = \sum_{i=0}^{p-1} t^i\, K_{M,i}(x) .
\end{equation}

\begin{definition}\label{def:int-curv}
Let $M$ be a compact oriented smooth hypersurface $M$ of $S^p$
and $U$ be an open subset of $M$.
The {\em integral $\mc_i(U)$ of $i$th curvature} and the
{\em integral $|\mc_i|(U)$ of $i$th absolute curvature} over $U$,
with respect to the ambient space~$M$, are defined as
($0\le i\le p-1$)
\begin{equation*}\label{eq:def-int-i-curv}
 \mc_i(U) := \int_U K_{M,i}\, dM,\qquad\mbox{and}\qquad
 |\mc_i|(U) := \int_U |K_{M,i}|\, dM .
\end{equation*}
\end{definition}

A few remarks:
$|\mu_i|$ is monotone in the sense that
$|\mc_i|(U_1)\leq |\mc_i|(U_2)$ for $U_1\subseteq U_2$.
%(which is not true for $\mu_i$).
Also, note that $|\mc_i(U)|\le |\mc_i|(U)$
and $\mc_i(\emptyset)=|\mc_i|(\emptyset)=0$.
Moreover,
$\mc_0(U)= |\mc_0(U)|=\vol_{p-1}(U)$.
By Example~\ref{ex:sphere},
$|\mc_i|(S^{p-1})=0$ for $i>0$.

\begin{example}\label{ex:geoball}
Consider the boundary $M_\a$ of the ball
$B_R(q,\a)$ in $S^p$ of radius $0<\a\le\pi/2$ centered at~$q$.
Clearly, $M_\a$ is an $(p-1)$-dimensional sphere of radius $\sin\a$
that is described by the equations
\begin{equation*}
 x_0 = \cos\a,\ x_1^2 +\cdots +x_p^2 = \sin^2\a
\end{equation*}
if $q=(1,0,\ldots,0)$. We orient $M_\a$ by the unit normal vector
field on $S^p$ pointing towards $q$. It is straight-forward to see
that the second fundamental form of $M_\a$ satisfies $\II_{M_\a}(x) =
(\cot\a)\, \mathrm{id}_{p-1}$ for all $x\in M_\a$.
%(cf. Figure~\ref{fig:Ma}).
%\begin{figure}
\begin{center}
 \input fig_real_2.pictex
\end{center}
%\caption{The manifold $M_\a$}\label{fig:Ma}
%\end{figure}
Hence all the principal curvatures of $M_\a$ at $x$ are equal to $\cot\a$.
Therefore the $i$th curvature of $M_\a$ satisfies
$K_{M_\a,i}(x) = {p-1\choose i}\,(\cot\a)^{i}$, a quantity independent of
$x\in M_\a$. For the integral of the $i$th curvature we obtain
\begin{equation}\label{eq:muiM}
 \mu_i(M_\a) = K_{M_\a,i}\, \vol_{p-1} M_\a = {p-1\choose i} \Oh_{p-1}
        (\sin\a)^{p-i-1} (\cos\a)^{i},
\end{equation}
using the fact that $\vol_{p-1} M_\a= \Oh_{p-1}(\sin\a)^{p-1}$.
In particular, note that
$\mu_{p-1}(M_\a) = \Oh_{p-1}(\cos\a)^{p-1}$.
Finally note that $\mu_i(U)=|\mu_i|(U)$ for all open subsets $U$ of $M_\a$,
since all the principal curvatures are nop-negative.
\end{example}

%%%%%%%%%%%%%%%%%%%%%%%%%%%%%%%%%%%%%%%%
\subsection{A kinematic formula from integral geometry for
spheres}\label{se:kin-sphere}

%We are thus left with the task of bounding $|\mu_i|(U)$.
%We will do so in \S\ref{sec:mu} below.

We recall here a basic formula of integral geometry.
The orthogonal group $G=O(p+1)$ operates on $S^p$ in the natural way.
We will denote by $dG$ the volume element on the compact Lie group $G$
normalized such that the volume of $G$ equals one.
We will interpret $S^{i}$ as a submanifold of $S^p$ for $i\le  p$,
e.g., given by the equations $x_{i+1}=\cdots=x_p=0$.

Let $M$ be a compact oriented smooth hypersurface of $S^p$. It
follows by standard methods from Sard's lemma~\cite{spiv1:79} that
$gM$ intersects $S^{i+1}$ transversally for almost all $g\in G$.
Hence, for almost all $g\in G$, the intersection $gM\cap S^{i+1}$ is
either empty or a smooth hypersurface of $S^{i+1}$. Moreover,
this intersection inherits an
orientation from~$M$ in a natural way as follows: let $\nml$ be the
distinguished unit normal vector field of~$M$. Then we require that
the distinguished unit normal vector of the hypersurface $gM\cap
S^{i+1}$ in $S^{i+1}$ at~$x$ lies in the positive halfspace of
$T_xM$ determined by $\nml$.

Therefore, for almost all $g\in G$, the integral of the $i$th
curvature $\mc_i(gM\cap S^{i+1})$
of $gM\cap S^{i+1}$, considered as submanifold of $S^{i+1}$,
is well defined, and this is also the case for $\mc_i(gU\cap S^{i+1})$
when $U$ denotes an open subset of $M$.

We will need the following special case of the principal kinematic
formula of integral geometry for spheres.
A proof can be found in~\cite{howa:93}.
For Euclidean space a corresponding result was stated by Chern~\cite{cher:66}.

\begin{theorem}\label{th:kin-form-S}
Let $U$ be an open subset of a compact oriented smooth hypersurface~$M$
of $S^p$ and $0\le i <p-1$. Then we have
\begin{equation*}
\mc_i(U) = \kc(p,i) \int_{g\in G} \mc_i(gU\cap S^{i+1})\, dG(g),
\end{equation*}
where
$
%\begin{equation}\epf
\kc(p,i) = (p-i-1){p-1\choose i}\,\frac{\Oh_{p-1}\Oh_{p}}
           {\Oh_{i}\Oh_{i+1}\Oh_{p-i-2}}.
%\end{equation}
$
\hfill$\Box$
\end{theorem}

While the existence of the constants $\kc(p,i)$ follows 
from~\cite{cher:66,howa:93}, it is quite cumbersome to extract 
explicit formulas for $\kc(p,i)$ from these sources.
This is partly due to the fact that the 
quantities $\mu_i$ (and even $\Oh_i$) have slightly
different meanings in the literature.
For the convenience of the reader, we have therefore included 
a short derivation of these constants in the appendix. 

For future reference, we state the case $i=0$ of 
Theorem~\ref{th:kin-form-S} in slightly more general form.

\begin{corollary}\label{cor:intgeo}
For any semialgebraic subset $T$ of $S^p$ such that $\dim T\le p-1$ we have
\begin{equation}%\epf
 \vol_{p-1} T = \frac{\Oh_{p-1}}{2} \int_{g\in G} |T \cap gS^{1}|\, dG(g).
\end{equation}
\end{corollary}

\begin{proof}
Using the comments given at the end of \S\ref{se:dist-vol-tub}, 
it is easy to reduce to the case where $T$ is an open subset of
a hypersurface of $S^p$. Now apply Theorem~\ref{th:kin-form-S}
for $i=0$, taking into account that neither the compactness 
nor the orientability assumption are 
necessary in that case, cf.~\cite[\S15.2]{sant:76}.
\end{proof}

%%%%%%%%%%%%%%%%%%%%%%%%%%%%%%%%%%%%
\section{On the volume of tubes around real algebraic sets}
\label{se:vol-tubes}

The goal of this section is to provide the proof of Theorem~\ref{th:realtubel}.

%%%%%%%%%%%%%%%%%%%%%%%%%%%%%%%%%%%%
\subsection{Bounding the volume of tubes of smooth hypersurfaces}
\label{se:weyl-etc}

In an important paper, Weyl~\cite{weyl:39} derived a formula for the
volume of tubes around a compact submanifold of Euclidean space or a
sphere. However, this formula only holds for a sufficiently small
radius. The following proposition gives an upper bound on the volume
of tubes around a hypersurface that holds for any radius. Compare
also Gray~\cite[Theorem~8.4, (8.6), p.~162]{gray:90}.

\begin{proposition}\label{pro:tube-vol-mc}
Let $M$ be a compact oriented smooth hypersurface of $S^p$ and
$U$ be an open subset of $M$.
Then we have for all $0<\a \le \pi/2$
\begin{equation*}
 \vol_p \Tusn(U,\a) \le 2\sum_{i=0}^{p-1} J_{p,i+1}(\a)\, |\mu_i|(U) .
\end{equation*}
%where $J_k(\a)$ is defined as
%\begin{equation*}
%  J_{k}(\a) := \int_0^\a (\sin\rho)^{k-1}\, (\cos\rho)^{p-k}\, d\rho .
%\end{equation*}
\end{proposition}

\begin{proof}
Let $\nml\colon M\to S^p$ be the unit normal vector field on $M$
corresponding to its orientation.
For $x\in M$ consider the following parametrization
\begin{equation*}
 \varphi_x\colon\R\to S^p,\  \varphi_x(t)
 = \frac{x + t\nml(x)}{\| x+ t\nml(x)\|}=\frac{x + t\nml(x)}{\left(1+t^2\right)^{\frac{1}{2}}}
\end{equation*}
of the half great circle intersecting $M$ at $x$ orthogonally.
Note that $d_R(x,\varphi_x(t))=\arctan t$.
%\begin{figure}
\begin{center}
  \input fig_real_3.pictex
\end{center}
%\caption{The point $\varphi_x(t)$ and the magnitudes $t$ and $\a$}\label{fig:paramet}
%\end{figure}
Consider the following differentiable map of Riemannian manifolds
\begin{equation*}
 \varphi\colon M\times\R \to S^p, (x,t)\mapsto\varphi_x(t).
\end{equation*}
Let $U$ be an open subset of $M$, $0<\a\le\pi/2$, and put $\t=\tan\a$.
We denote by $\Tusplus(U,\a)$ and $\Tusminus(U,\a)$ the images of $U\times (0,\t)$
and $U\times (-\t,0)$ under the map~$\varphi$, respectively.
Clearly, $\Tusn(U,\a)=U\cup\Tusplus(U,\a)\cup\Tusminus(U,\a)$.

We apply the transformation formula of Proposition~\ref{pro:trafo}
to the surjective differentiable map
$\varphi\colon U\times (0,\tau) \to \Tusplus(U,\a)$
of Riemannian manifolds. This yields
\begin{equation*}
 \int_{U\times (0,\t)} |\det D\varphi|\, d(M\times\R)
 = \int_{y\in\Tusplus(U,\a)} |\varphi^{-1}(y)|\, dS^p
   \ \ge\  \vol_p \Tusplus(U,\a).
\end{equation*}
By Fubini,
$\int_{U\times (0,\t)} |\det D\varphi| d(M\times\R) = \int_0^\t g(t)\,dt$, where
$$
 g(t) :=\int_{x\in U} |\det D\varphi|(x,t) \, dM(x) .
$$

\medskip
\noindent{\bf Claim~A.} The determinant of the derivative
$D\varphi(x,t)$ of $\varphi$ at $(x,t)\in M\times\R$ satisfies
\begin{equation}
 |\det D\varphi(x,t)| = \frac1{(1+t^2)^{(p+1)/2}}\;
 \bigl|\det(\mathrm{id}_{T_xM} - t L_M(x))\bigr| .
\end{equation}
Using this claim, whose proof is postponed to the end,
we obtain
\begin{eqnarray*}
g(t)  &=& \int_{x\in U} \frac{1}{(1+t^2)^{\frac{p+1}{2}}}
              |\det(\mathrm{id}_{T_xM} - tL_M(x))|\, dM(x)
              \qquad\mbox{(by Claim~A)}\\
    &\leq& \sum_{i=0}^{p-1} \frac{|t|^i}{(1+t^2)^{\frac{p+1}{2}}}
           \int_U |K_{M,i}|\,dM \qquad\mbox{(by (\ref{eq:charpol}))}\\
    &=& \sum_{i=0}^{p-1} \frac{|t|^i}{(1+t^2)^{\frac{p+1}{2}}}\ 
             |\mu_i|(U).
\end{eqnarray*}
By making the substitution $t=\tan \rho$ (recall $\t=\tan\a$) 
we get
\begin{equation*}
 \int_0^\t \frac{t^i}{(1+t^2)^{(p+1)/2}}\, dt =
 \int_0^\a (\cos\rho)^{p-i-1} (\sin\rho)^{i}\,d\rho
 = J_{p,i+1}(\a).
\end{equation*}
Altogether we conclude that
$$
\vol_p \Tusplus(U,\a) \le  \int_{0}^\t g(t)\, dt
\leq \sum_{i=0}^{p-1} J_{p,i+1}(\alpha)|\mu_i|(U).
$$
The same estimate can be shown for $\vol_p \Tusminus(U,\a)$, which
implies the desired estimate  of $\vol_p \Tusn(U,\a)$.

It remains to prove Claim~A.
Choose a local parametrization
$x=x(v)\in\R^{p+1}$ of~$M$ in terms of coordinates $v_1,\ldots,v_{p-1}$
defined in a neighborhood of $0$. We assume that
$\partial_{v_1} x,\ldots,\partial_{v_{p-1}} x$ are orthonormal at $0$.
Abusing notation we will interpret $\nml=\nml(v)$ as a function of $v$.
The matrix $(\lambda_{ij})$ of $L_M$ with respect to the basis
$\partial_{v_j}x$ of $T_xM$ is given by
$-\partial_{v_i}\nml = \sum_j \lambda_{ij}\partial_{v_j}x$.

Consider the map $(v,t)\mapsto R(v,t):= x(v) + t\nml(v)\in\R^{p+1}$.
Then
\begin{equation*}
(v,t)\mapsto \psi(v,t):=\varphi(x(v),t)=
\frac{R(v,t)}{(1+t^2)^{\frac12}}
\end{equation*}
provides a local parametrization of $S^p$. In the following let $[R,
\partial_{t}R, \partial_{v_1}R,\ldots,
\partial_{v_{p-1}}R ]$ denote the square matrix of size $p+1$ whose
rows are $R$ and the partial derivatives of~$R$. Using the
%fact that $\|\psi(v,t)\|=1$ and $\<\psi,\partial_t\psi\>=\<\psi,\partial_{v_j}\psi\>=0$
%for $1\leq j\leq p-1$, it is easy to check that
multilinearity of the determinant and the fact
$$
 \partial_t\psi = (1+t^2)^{-1/2}\partial_t R - t (1+t^2)^{-3/2}R
$$
we obtain by a short calculation that
\begin{align*}
 |\det D\psi(v,t)| &= |\det
 [\psi,\partial_t\psi,\partial_{v_1}\psi,\ldots,\partial_{v_{p-1}}\psi]|\\
  &=\frac1{(1+t^2)^{(p+1)/2}}\;
 \bigl|\det [R, \partial_{t}R, \partial_{v_1}R,\ldots, \partial_{v_{p-1}}R
  ]\bigr|.
\end{align*}
Computing partial derivatives we get
$ \partial_{v_i}R = \partial_{v_i}x + t \partial_{v_i} \nml$ and $\partial_{t}R = \nml$.
Using $\partial_{v_i}\nml = - \sum_j \lambda_{ij}\partial_{v_j}x$,
we get from this 
\begin{equation*}
\mbox{$\partial_{v_i}R = \sum_j (\delta_{ij} -t \lambda_{ij})\partial_{v_j}x$.}
\end{equation*}
Hence we obtain
\begin{eqnarray*}
  \det [R, \partial_{t}R, \partial_{v_1}R,\ldots, \partial_{v_{p-1}}R ]
 &=& \det [x + t\nml, \nml,
     \partial_{v_1}R,\ldots, \partial_{v_{p-1}}R ] \\
 &=& \det [x, \nml, \partial_{v_1}R,\ldots, \partial_{v_{p-1}}R ] \\
 &=& \det(\delta_{ij} -t \lambda_{ij})
    \det [x, \nml, \partial_{v_1}x,\ldots, \partial_{v_{p-1}}x ] .
\end{eqnarray*}
Since we assume that $\partial_{v_1} x,\ldots,\partial_{v_{p-1}} x$
are orthonormal at $v=0$, we conclude (using the chain rule
$D\psi=D\varphi Dx$) that
\begin{equation*}
  |\det D\varphi(x,t)| =  |\det D\psi(0,t)|
  = \frac1{(1+t^2)^{(p+1)/2}} |\det(\delta_{ij} -t \lambda_{ij})|,
\end{equation*}
which shows Claim~A.
\end{proof}

%%%%%%%%%%%%%%%%%%%%%%%%%%%%%%%%%%%%%%%%%%%%%%%%%%%%%%%
\subsection{Bounding integrals of absolute curvature in terms of degree}
\label{sec:mu}

In this section let $f\in\R[X_0,\ldots,X_p]$ be homogeneous of degree $d>0$
with nonempty zero set $V\subseteq S^p$ such that the derivative of
the restriction of $f$ to $S^p$ does not vanish on $V$.
%$f\colon S^p\to\R$ does not vanish on~$V$.
Then $V$ is a compact smooth hypersurface of $S^p$.
We orient $V$ by the following unit normal vector field
(Gauss map)
\begin{equation*}
   \nml\colon V\to S^p,\
   \nml(x) = \|\grad f(x)\|^{-1}\,\grad f(x) .
\end{equation*}
The goal of this section is to bound the integrals of absolute curvature
on patches of~$V$.

\begin{proposition}\label{le:mc-estimate}
For $a\in S^p$, $0<\s\le 1$, and $0\le i < p$ we have
\begin{equation*}
 |\mc_i|(V\cap B_\P(a,\s))\ \le\
 2{p-1\choose i} \Oh_{p-1}\, d^{i+1}\, \s^{p-i-1} .
\end{equation*}
\end{proposition}

The proof is based on the following lemma.
%We first show how to derive this statement from the following result.
%The next result is the first step towards bounding integrals
%of curvatures in terms of the degree.

\begin{lemma}\label{le:abs-gauss-curv2}
We have $|\mu_{p-1}|(V)\le \Oh_{p-1} d^{p}$.
\end{lemma}

\begin{proof}
Recall that the determinant of the linear map
$L_V(x)=-D\nml(x)\colon T_x V\to T_x V$
equals $K_{V,p-1}$ (cf.~\S\ref{se:dg-sphere}).
We may assume without loss of generality that the
open subset $U:=\{x\in V\mid \rank(D\nml(x)) =p-1\}$ of $V$ is
nonempty (otherwise $\mu_{p-1}(V)=0$).
We would like to apply Proposition~\ref{pro:trafo}
to the restriction of $\nu$ to~$U$, but face the problem that
it is only an immersion so that $\nu(U)$ might not be a submanifold of $S^p$.
In order to circumvent this, we use some standard facts
of real algebraic geometry~\cite{bocr:98}.

Consider the Zariski closure $W$ of $\nu(U)$ in $S^p$,
which is a real algebraic variety of dimension $p-1$.
Its regular locus $W_1$ is a submanifold of $S^p$ of dimension $p-1$.
Consider the open subset $V_1 := U\cap \nu^{-1}(W_1)$ of $V$
and the restriction $\nu_1\colon V_1\to W_1$ of~$\nu$.
The singular locus $\mathrm{Sing}(W)=W\setminus W_1$ is an algebraic subset
of dimension strictly less than $\dim W$.
Since $\nu_1$ is an immersion, we conclude that
$U\setminus V_1 = \nu_1^{-1}(\mathrm{Sing}(W))$
has dimension strictly less than $p-1$.
We therefore obtain
$$
 |\mu_{p-1}(V)| = |\mu_{p-1}(U)| = |\mu_{p-1}(V_1)|.
$$
Applying the transformation formula of Proposition~\ref{pro:trafo}
to $\nu_1$ we get
$$
 |\mu_{p-1}(V)| = \int_{V_1} |\det D\nu_1|\, dV
  = \int_{y\in W_1} |\nu_1^{-1}(y)|\, dW_1(y) .
$$
Consider for $\ell\in\N\cup\{\infty\}$ the semialgebraic sets
$F_\ell:=\{y\in W_1\mid |\nu_1^{-1}(y)|=\ell \}$.
Since the above integral is finite, $\vol_{p-1}F_\infty=0$, and
therefore $\dim F_\infty < p-1$.
We obtain
$$
\int_{y\in W_1} |\nu_1^{-1}(y)|\, dW_1(y)
 = \sum_{\ell\ge 0} \ell\,\vol_{p-1}\,F_\ell .
$$
Corollary~\ref{cor:intgeo} applied to $F_\ell$ yields the following
\begin{equation*}
 \vol_{p-1} F_\ell = \frac{\Oh_{p-1}}{2}\
 \int_{g\in G}\, |F_\ell \cap gS^1|\, dG(g).
\end{equation*}
Combining these findings we get
\begin{equation*}\label{eq:Knest}
 |\mu_{p-1}|(V)=
 \frac{\Oh_{p-1}}{2}\ \int_{g\in G}\,
 \sum_{\ell\ge 0} \ell\,|F_\ell \cap gS^1|\,dG(g)
 = \frac{\Oh_{p-1}}{2}\ \int_{g\in G}\, |\nu_1^{-1}(gS^1)| \,dG(g).
\end{equation*}
(In order to see the last equality use that $gS^1$ does not
intersect $F_\infty$ almost surely.)

It is now sufficient to prove that
$|\nu^{-1}(gS^1)| \le 2d^p$
for almost all $g\in G$.
To simplify notation suppose first that $g=\mathrm{id}$.
%After applying a linear coordinate transformation to~$f$
%we may assume without loss of generality that
Let $y_0,\ldots,y_p$ denote coordinate functions on $\R^{p+1}$ and
$S^1=\{y\in S^p \mid  y_2=\cdots =y_p=0\}$.
A point $x\in\R^{p+1}$ lies in $\nu^{-1}(S^1)$ iff it satisfies
the following system of equations
\begin{equation*}
 \sum_i x_i^2 -1 =0,\ f(x)=0,\ \partial_2 f(x) =\cdots=\partial_p f(x)= 0.
\end{equation*}
If all real solutions of this system of equations are nondegenerated,
then they are isolated in $\C^{p+1}$.
By B\'ezout's theorem~\cite{shaf:74} the number of these solutions
is bounded by $2\,d\,(d-1)^{p-1}\le 2d^p$.
However, one can show along the lines in~\cite{miln:64}, that
the nondegeneracy condition is satisfied for almost
all $g\in G$.
This finishes the proof.
\end{proof}

\proofof{Proposition \ref{le:mc-estimate}}
Put $U:=V\cap B_\P(a,\s)$. The case $i=p-1$ is already settled by
Lemma~\ref{le:abs-gauss-curv2}, as $|\mu_i|(U)\le |\mu_i|(V)$.
So we may assume $i<p-1$.

Let $U_+$ be the set of points of $U$
where $K_{V,i}$ is positive and similarly define $U_-$ where
$K_{V,i}$ is negative. Then $|\mu_i|(U) = |\mu_i(U_+)| +|\mu_i(U_-)|$.

Let $g\in G=O(p+1)$ such that $V$ intersects $gS^{i+1}$ transversally
and such that the intersection is nonempty.
Then $V\cap gS^{i+1}$ is the zero set in $gS^{i+1}$ of
the homogeneous polynomial~$f$ of degree~$d$.
By transversality, the derivative of the restriction of $f$ to
$gS^{i+1}$ does not vanish on $V\cap gS^{i+1}$.
Hence we may apply Lemma~\ref{le:abs-gauss-curv2},
which gives the estimate
$$
 |\mu_i|(V\cap gS^{i+1}) \leq \Oh_i d^{i+1}.
$$
By the monotonicity of $|\mu_i|$
(now refering to the hypersurface $V\cap gS^{i+1}$ of $gS^{i+1}$)
we have
\begin{equation*}
  |\mu_i(U_+\cap gS^{i+1})| \le
  |\mu_i|(U_+\cap gS^{i+1}) \le |\mu_i|(V\cap gS^{i+1})
\leq \Oh_i d^{i+1}.
\end{equation*}
The Kinematic Formula of Theorem~\ref{th:kin-form-S} implies that
$$
  |\mc_i(U_+)| 
  \le \kc(p,i) \int_{g\in G} |\mc_i(gU_+\cap S^{i+1})| \, dG(g).
$$
Therefore, we obtain
\begin{equation*}
 |\mc_i(U_+)|\ \le\ \kc(p,i)\ \Oh_i d^{i+1}\
  \Prob_{g\in G} \{gU_+\cap S^{i+1}\ne\emptyset\},
\end{equation*}
where the probability is taken with respect to the
uniform distribution in $G$.
% and similarly for $U_-$.
We may estimate this probability as follows
(put $\a=\arcsin\s$)
\begin{eqnarray*}
\lefteqn{\Prob_{g\in G} \{B_\P(ga,\s)\cap S^{i+1}\ne\emptyset\}
  =  \Prob_{a'\in S^p} \{B_\P(a',\s)\cap S^{i+1}\ne\emptyset\} }\\
 &=& \frac{1}{\Oh_{p}}\vol_p T_\P(S^{i+1},\s)
 = \frac{\Oh_{i+1}\Oh_{p-i-2}}{\Oh_{p}}\, J_{p,p-i-1}(\a)
 \le\frac{\Oh_{i+1}\Oh_{p-i-2}}{\Oh_{p}}\ \frac{\s^{p-i-1}}{p-i-1},
\end{eqnarray*}
where we used Lemmas~\ref{le:J1} and~\ref{le:J2}
for the last two steps.
Multiplying this with the formula for $\kc(p,i)$,
the expression miraculously simplifies and we get
$$
 |\mc_i|(U_+)\ \le\ {p-1\choose i} \Oh_{p-1}\ d^{i+1}\ \s^{p-i-1} .
$$
The same estimate can be shown for $ |\mc_i|(U_-)$, which
proves the proposition.
\endproofof

%%%%%%%%%%%%%%%%%%%%%%%%%%%%%%%%%%%%%%%%%%
\subsection{Extension to the non-smooth case: proof of 
Theorem~\ref{th:realtubel}}\label{se:realtubelproof}

The following proposition estimates the volume of the tube around a
patch of a smooth hypersurface in the sphere.

\begin{proposition}\label{pro:realtubel}
Let $f\in\R[X_0,\ldots,X_p]$ be homogeneous of even degree $d>0$
with zero set $V=\mZ(f)$ in $S^p$.
Assume that the derivative of the restriction of $f$ to $S^p$
does not vanish on $V$.
(Thus $V$ is a smooth hypersurface in $S^p$.)
Let $a\in S^p$ and $0<\e,\s\le 1$. Then
\begin{equation*}
 \vol_p\Tuproj(V\cap B_{\P}(a,\s),\e) \ \le\
 \frac{4\Oh_{p-1}}{p} \sum_{k=1}^{p-1} {p\choose k}\, d^{k}\,
  \e^{k}\,\sigma^{p-k}+2\Oh_{p}\, d^p\,\e^p.
\end{equation*}
\end{proposition}

\begin{proof}
Put $U:=V\cap B_\P(a,\s)$.
Take $0<\a \le\pi/2$ such that $\e=\sin\a$.
Proposition~\ref{pro:tube-vol-mc} implies
\begin{equation*}
 \vol_p \Tuproj(U,\e) = \vol_p \Tusn(U,\a)
 \le 2\sum_{i=0}^{p-1} J_{p,i+1}(\a)\, |\mu_i|(U) .
\end{equation*}
Combining this with Proposition~\ref{le:mc-estimate}
we obtain
\begin{equation*}
  \vol_p\Tuproj(U,\e) \le\
 4\sum_{i=0}^{p-1} {p-1\choose i} \Oh_{p-1} d^{i+1} \sigma^{p-i-1}
 J_{p,i+1}(\a).
\end{equation*}
Using the estimates of Lemma~\ref{le:J2} we obtain (put $k=i+1$ and consider
separately the term for $k=p$)
\begin{equation*}
 \vol_p\Tuproj(U,\e) \le\
 4 \sum_{k=1}^{p-1} {p-1\choose k-1} \Oh_{p-1}\, d^k\,
  \sigma^{p-k}\, \frac{\e^k}{k}
  + 4\Oh_{p-1}\, d^p \frac{\Oh_p}{2\Oh_{p-1}}\, \e^p.
\end{equation*}
Now use ${p-1\choose k-1}=\frac{k}{p}{p\choose k}$ to get the
desired upper bound on $\vol_p \Tuproj(U,\e)$.
\end{proof}

\proofof{Theorem~\ref{th:realtubel}}
We have to remove the smoothness assumption in
Proposition~\ref{pro:realtubel} and to estimate the volume
of the $\e$-neighborhood instead of the $\e$-tube.

Assume $W=\mZ(f_1,\ldots,f_r)$ with homogeneous polynomials $f_i$
of degree~$d_i$. Then $W$ is the zero set in $S^p$ of the polynomial
\begin{equation*}
   f(X):=\sum_{i=1}^r f_i(X)^2\|X\|^{2d-2d_i},
\end{equation*}
which is homogeneous of degree~$2d$.
Our assumption $W\ne S^p$ implies $\dim W<p$.

Let $\d>0$ be smaller than any positive critical value of the
restriction $\tilde{f}\colon S^p\to\R$ of $f$ to $S^p$. Then
$D_\d:=\{\xi\in S^p\mid \tilde{f}(\xi)\le \d\}$
is a compact domain in $S^p$ with smooth boundary
\begin{equation*}
   \partial D_\d =\{\xi\in S^p\mid \tilde{f}(\xi)=\d \}.
\end{equation*}
Indeed, the derivative of $\tilde{f}-\d$
does not vanish on $\partial D_\d$
(use $\sum_i x_i\partial_i f(x) = 2d f(x)$).
Moreover, note that $W=\cap_{\d>0} D_\d$ and hence
$\lim_{\d\to 0}\vol_p D_\d = \vol_p(W) =0$, since
$\dim W<p$.

\bigskip

\noindent{\bf Claim B.} We have $T_{\P}(W,\e)\subseteq D_\d\cup
T_{\P}(\partial D_\d,\e)$ for $0<\e\le 1$.
\medskip

In order to see this, let $x\in T_{\P}(W,\e)\setminus D_\d$ and
$\g\colon [0,1]\to S^p$ be a segment of Riemannian length less
than $\arcsin\e$ such that $\g(1)=x$ and $\g(0)\in W$. Consider
$F\colon[0,1]\to\R, F(t):= \tilde{f}(\gamma(t))$. By assumption
$F(1)=\tilde{f}(x)>\d$ and $F(0)=0$. Hence there exists $\tau\in (0,1)$ such that
$F(\tau)=\d$. Thus $\g(\tau)\in\partial D_\d$ and
$d_{\P}(x,\partial D_\d)\le d_{\P}(x,\gamma(\tau))<\e$, which shows the claim.
\medskip

Using the triangle inequality for the projective distance,
it is easy to see that (cf.~Figure~4)
\begin{equation}\label{eq:inclusion}
 T_{\P}(\partial D_\d,\e)\cap B_{\P}(a,\s) \subseteq
 \Tuproj(\partial D_\d\cap B_{\P}(a,\s+\e),\e).
\end{equation}
\begin{center}
   \input fig_real_5.pictex
\end{center}
Combining (\ref{eq:inclusion}) with Claim~B we arrive at
\begin{equation*}
 T_{\P}(W,\e)\cap B_{\P}(a,\s) \subseteq D_\d \cup
 \Tuproj(\partial D_\d\cap B_{\P}(a,\s+\e),\e) .
\end{equation*}
We apply Proposition~\ref{pro:realtubel} to $V=\partial D_\d
=\mZ(f-\d\|x\|^{2d})$ intersected with the ball $B_{\P}(a,\s+\e)$.
This implies
$$
 \vol_p \Tuproj(\partial D_\d\cap B_{\P}(a,\s+\e),\e)
 \le \frac{4\Oh_{p-1}}{p} \sum_{k=1}^{p-1} {p\choose k}\, (2d)^{k}\,
  \e^{k}\,(\sigma + \e)^{p-k} + 2\Oh_{p}\, (2d)^p\,(\sigma + \e)^p.
$$
Taking into account that $\vol_p B_{\P}(a,\s) \ge
\Oh_{p-1}\frac{\s^p}{p}$ (cf.~Lemmas~\ref{le:J2} and~\ref{le:J1})
we obtain
\begin{align*}
 \frac{\vol_p \left(T_{\P}(W,\e)\cap B_{\P}(a,\s)\right)}
   {\vol_p B_{\P}(a,\s)}\ & \le\
   \frac{\vol_p D_\d}{\vol_p B_{\P}(a,\s)}\
   + \ \frac{\vol_p\Tuproj(\partial D_\d\cap B_{\P}(a,\s+\e),\e)}
     {\vol_p B_{\P}(a,\s)} \\
   &\le\
  \frac{\vol_p D_\d}{\vol_p B_{\P}(a,\s)} \ +\
   4\sum_{k=1}^{p-1} {p\choose k}\, (2d)^{k}\, \left(1
   +\frac{\e}{\s}\right)^{p-k}\left(\frac{\e}{\s}\right)^{k}\\
   & \quad + \frac{2p\Oh_p}{\Oh_{p-1}}\, (2d)^p\,\left(\frac{\e}{\s}\right)^{p}.
\end{align*}
Taking the limit for $\d\to 0$ the first term vanishes and
the assertion follows.
\proofend

%%%%%%%%%%%%%%%%%%%%%%%%%%%%%%%%%%%%%%%%%%%%%%%%%%%%%%%%%%%%%
\subsection{Estimating expected values: proof of Theorem~\ref{th:main}}

The tail bounds in Theorem~\ref{th:main} follow from
Theorem~\ref{th:realtubel} as indicated in the introduction.
It thus suffices to show how to derive the claim on expected values
from the tail bounds.
This is achieved by the following proposition.

We use the inequality
$\frac{\Oh_p}{2\Oh_{p-1}} -\frac1{p} \leq \frac12$,
valid for $p\ge 2$, which implies
$\frac{p\Oh_p}{\Oh_{p-1}} \leq p+2$.

\begin{proposition}
For $0<\sigma\leq 1$ let $X_\sigma\geq 1$ be a random variable
satisfying, for all $0<\e \leq 1$ and $p\ge 2$,
$$
   \Prob\{X_\sigma\geq 1/\e\}\leq
   4\sum_{k=1}^{p-1} {p\choose k} (2d)^k\,
   \left(1 +\frac{\e}{\s}\right)^{p-k} \left(\frac{\e}{\s}\right)^{k}
   +2(p+2)\, (2d)^p\,\left(\frac{\e}{\s}\right)^{p}.
$$
Then, for $\e\leq \frac{\sigma}{(1+2d)(p-1)}$, we have
$\Prob\{X_\sigma\geq 1/\e\}\leq
(8e+4)dp \frac{\e}{\sigma}$ and
$$
   \bfE(\ln X_\sigma)\leq   2\ln p+2\ln d+2\ln\frac{1}{\sigma}+5.5.
$$
\end{proposition}

\begin{Proof}
$\Prob\{X_\sigma\geq 1/\e\}$ is bounded by
\begin{align*}
%   \Prob\{&X_\sigma\geq 1/\e\}\\
% &\leq 4\left[\sum_{k=1}^{p-1} {p\choose k} (2d)^k\,
%   \left(1 +\frac{\e}{\s}\right)^{p-k}\left(\frac{\e}{\s}\right)^{k}
%   + (\frac{p}{2}+1)\, (2d)^p\,\left(\frac{\e}{\s}\right)^{p}\right]\\
 &4\left[\sum_{k=1}^{p} {p\choose k} (2d)^k\,
   \left(1 +\frac{\e}{\s}\right)^{p-k}\left(\frac{\e}{\s}\right)^{k}
   + \frac{p}{2}\,  (2d)^p\,\left(\frac{\e}{\s}\right)^{p}\right]\\
 &=\; 8d\frac{\e}{\sigma}\left[
   \sum_{k=1}^{p} {p\choose k} (2d)^{k-1}\,
   \left(1 +\frac{\e}{\s}\right)^{p-k}\left(\frac{\e}{\s}\right)^{k-1}
   + \frac{p}{2} \, (2d)^{p-1} \left(\frac{\e}{\s}\right)^{p-1}\right]\\
 &\leq\; \frac{8dp\e}{\sigma}\left[
   \sum_{k=0}^{p-1} {p-1\choose k} (2d)^k\,
   \left(1 +\frac{\e}{\s}\right)^{p-1-k}
   \left(\frac{\e}{\s}\right)^{k}
   +\frac12 (2d)^{p-1}\,\left(\frac{\e}{\s}\right)^{p-1}\right]
\end{align*}
Using that $\e\leq \frac{\sigma}{(1+2d)(p-1)}$, this can be further bounded
to obtain
%\begin{align*}
$$
 \Prob\{X_\sigma\geq 1/\e\} \le
 \frac{8dp\e}{\sigma}\left[ \left(1 +\frac{1}{p-1}\right)^{p-1} 
  +\frac12 \right]
\leq\; \frac{8dp}{\sigma}(e+\frac12)\, \e.
%\qquad\mbox{since $\e\leq \frac{\sigma}{(1+2d)(p-1)}$}\\
%\end{align*}
$$
For the bound on the expectation we use
Proposition~2.4 in~\cite{BCL:06a} which implies 
\begin{eqnarray*}
  \bfE(\ln X_\sigma)&\leq& \ln\left(\frac{(1+2d)(p-1)}{\s}\right)
   +\ln\left(\frac{8dp}{\sigma}(e+\frac12)\right)+1 \\
 &\leq&
  2\ln p+2\ln d+2\ln\frac{1}{\sigma}+\ln(3\cdot 8(e+\frac12)e)
  \qquad\mbox{since $1+2d\leq 3d$}\\
 &\leq&
  2\ln p+2\ln d+2\ln\frac{1}{\sigma}+5.5,
\end{eqnarray*}
which completes the proof.
\end{Proof}

%%%%%%%%%%%%%%%%%%%%%%%%%%%%%%%%%%%%%%%%%%%%%%%%%%%%%%%%%%%%%%%%%
\section{Applications}\label{sec:appl}

We give several applications of Theorem~\ref{th:main}
to smooth analysis estimates for the condition numbers of
the following problems:
linear equation solving,
eigenvalue computations,
real polynomial equation solving, 
and zero counting.

%%%%%%%%%%%%%%%%%%%%%%%%%%%%%%
\subsection{Linear equation solving}\label{ss:les}

The first natural application of our result is for
the classical condition number $\kappa(A)$.
In~\cite{Wsch:04}, M.~Wschebor showed (solving
a conjecture posed in~\cite{ST:02}) that, for
all $n\times n$ real matrices $A$ with $\|A\|\leq 1$,
all $0<\sigma\leq 1$ and all $t>0$
$$
   \Prob_{Z\in N^{n^2}(A,\sigma^2)}(\kappa(Z)\geq t)\leq
   \frac{Kn}{\sigma t}
$$
with $K$ a universal constant.
Hereby, $\|A\|$  stands for the operator norm with respect to euclidean norm.
Note that, by Proposition~2.4 in~\cite{BCL:06a}, this implies
$$
   \sup_{\|A\|\le 1 }\ \bfE_{Z\in N^{n^2}(A,\sigma^2)}(\ln \kappa(Z))\leq
       \ln n+\ln\frac1{\sigma}+\ln K +1.
$$
We next compare Wschebor's result with
what can be obtained from Theorem~\ref{th:main}.
Let $\|A \|_F$ denote the Frobenius norm of a matrix $A\in\R^{n\times n}$,
which is induced by the inner product
$(A,B)\mapsto{\mathsf{trace}}(A B^T)$.
We have
$$
  \kappa(A)=\|A\|\|A^{-1}\|\leq \|A\|_F\|A^{-1}\|=:\kappa_F(A) .
$$
The Condition Number Theorem of Eckart and Young~\cite{EckYou}
states that
$\|A^{-1}\|=d_F(A,\Sigma)^{-1}$
where $\Sigma\subseteq\R^{n\times n}$ denotes the set of singular matrices
and $d_F$ is the distance induced by the Frobenius norm
(see also~\cite[Thm.~1, Ch.~11]{bcss:95}).
It follows that $\kappa_F$ is a conic condition number.
We can thus give upper bounds for
$\kappa_F(A)$ and they will hold as well for $\kappa(A)$.

\begin{corollary}\label{prop:cuad}
For all $n\geq 1$, $0<\sigma\le 1$, we have
$$
   \sup_{Z\in S^{n^2-1}}\bfE_{Z\in B_{\P}(A,\sigma)}(\ln \kappa_F(Z))\leq
    6\ln n+2\ln \frac1\sigma +5.5,
$$
where the expectation is over all $Z$ uniformly distributed
in the disk of radius $\sigma$ centered at $A$ in the sphere $S^{n^2-1}$
endowed with the projective distance.
\end{corollary}

\begin{Proof}
The variety $\Sigma$ of singular matrices is the zero set of the determinant,
which is a homogeneous polynomial of degree~$n$.
We now apply Theorem~\ref{th:main} with $p=n^2-1$.
\end{Proof}

Note, the bound in Corollary~\ref{prop:cuad} is of the same
order of magnitude than Wschebor's, worse by just a constant
factor. On the other hand, its derivation from
Theorem~\ref{th:main} is rather immediate.
We next extend this bound to rectangular matrices.

%%%%%%%%%%%%%%%%%%%%%%%%%%%%%%%%%%%%%%%%%%%%%%%%%%%%
\subsection{Moore-Penrose inversion}\label{ss:MPi}

Let $\ell\geq m$ and consider the space $\R^{\ell\times m}$ of
$\ell\times m$ rectangular matrices. Denote by
$\Sigma\subset\R^{\ell\x m}$ the subset of rank-deficient matrices.
For $A\not\in\Sigma$ let $A^\dagger$ denote its Moore-Penrose
inverse (see, e.g.,~\cite{BG,CM}). The condition number of $A$ (for
the computation of $A^\dagger$) is defined as
$$
  \cond^{\dagger}(A)=\lim_{\e\rightarrow 0}
  \sup_{\|\Delta A\|\leq \e}
  \frac{\|(A+\Delta A)^{\dagger} - A^{\dagger}\|\cdot \|A\|}
       {\|A^{\dagger}\|\cdot \|\Delta A\|}.
$$
This is not a conic condition number but it happens to be close to one.
One defines $\kappa^{\dagger}(A)=\|A\| \|A^\dagger\|$
and, since $\|A^\dagger\|=d_F(A,\Sigma)^{-1}$~\cite{GoLoan}, we
obtain
$$
  \kappa^{\dagger}(A)=\frac{\|A\|}{d_F(A,\Sigma)}.
$$
In addition (see~\cite[\S III.3]{Stewart}),
$$
  \kappa^{\dagger}(A)\leq\cond^{\dagger}(A)\leq
  \frac{1+\sqrt{5}}{2}\kappa^{\dagger}(A).
$$
Thus, $\ln(\cond^{\dagger}(A))$ differs from
$\ln(\kappa^{\dagger}(A))$
just by a small additive constant. As for square matrices,
$\kappa^{\dagger}(A)$ is not conic since the operator
norm is not induced by an inner product in
$\R^{\ell\times m}$. But, again, we can bound
$\kappa^{\dagger}(A)$ by the conic
condition number
$\kappa_F^{\dagger}(A):=\|A\|_F\|A^{\dagger}\|\ge\kappa^\dagger(A)$.

\begin{corollary}\label{prop:MP}
For all $\ell\geq m\geq 1$ and $0<\sigma\le 1$  we have
$$
   \sup_{A\in S^{\ell m-1}}
   \bfE_{Z\in B_{\P}(A,\sigma)}
   (\ln\kappa_{F}^{\dagger}(Z))\leq
   2\ln \ell + 4\ln m
   +2\ln\frac1\sigma + 5.5.
$$
\end{corollary}

\begin{Proof}
If a matrix $A$ is rank deficient then $\det A_0=0$
where $A_0$ is the $m\times m$ matrix obtained by removing
all rows of $A$ with index greater than $m$. Therefore
$\Sigma\subseteq\Sigma_0=\{A\in\R^{\ell\times m}\mid
\det A_0=0\}$. This implies that
$ \kappa_F^{\dagger}(A)\leq \frac{1}{d_F(A,\Sigma_0)}$
for $\|A\|_F=1$.
Since $\Sigma_0$ is the zero set of a homogeneous polynomial
of degree~$m$, an immediate application of
Theorem~\ref{th:main} with $W=\Sigma_0$ yields the claimed bound.
\end{Proof}

%%%%%%%%%%%%%%%%%%%%%%%%%%%%%%%%%%%%%%%%%%%%%%%%%%
\subsection{Real Eigenvalue computations}\label{ss:eigen}

Let $A\in\R^{n\times n}$ and $\lambda\in\R$ be a simple real 
eigenvalue of $A$.
Suppose that $x\in\R^n$ and $y\in\R^n$ are right and left 
eigenvectors associated
to $\lambda$, respectively (i.e., nonzero and satisfying
$Ax=\lambda x$ and $y^TA=\lambda y^T$).
>From the fact that $\lambda$ is a simple eigenvalue
one can deduce that $\langle x,y\rangle \ne 0$,
cf.\ Wilkinson~\cite{Wilkinson72}.

For any sufficiently small perturbation $\Delta A\in\R^{n\times n}$
there exists a unique real eigenvalue $\lambda+\Delta\lambda$
of $A+\Delta A$ close to $\lambda$. We thus have
$$
  (A+\Delta A)(x+\Delta x) = (\lambda + \Delta\lambda)(x+\Delta x),
$$
which implies up to second order terms
$\Delta A\, x+  A\, \Delta x \approx  \Delta \lambda\, x 
+  \lambda\, \Delta x$.
By multiplying with $y^T$ from the left we get
$$
 \Delta\lambda = \frac1{\langle x,y \rangle}y^T \Delta A x 
  + o(\| \Delta A\|).
$$
Moreover, $\sup_{\|\Delta A\|_F\le 1} |y^T\Delta A x|=\|x\|\,\|y\|$.

It therefore makes sense to define the condition number of $A$
for the computation of $\lambda$ as follows
\begin{equation}\label{eq:cnlambda}
   \kappa(A,\lambda) := \frac{ \|x\|\, \|y\|}{|\langle x,y \rangle|}
\end{equation}
and to set  $\kappa(A,\lambda) := \infty$ if $\lambda$ is a multiple 
eigenvalue of~$A$. 

If $A$ has real eigenvalues, we could define the 
condition number of $A$ for real eigenvalue computation as the maximum 
of $\kappa(A,\lambda)$ 
over all the real eigenvalues $\lambda$ of $A$.
In order to arrive at a reasonable definition that makes sense 
also for matrices having no real eigenvalues, consider 
the subset $\Sigma\subseteq\R^{n\times n}$ of matrices having 
a real multiple eigenvalue as the set of ill-posed inputs.  

We define the condition number of $A\in\R^{n\times n}$ 
for real eigenvalue computation by
\begin{equation}\label{eq:k-ev}
   \kappa_{{\sf eigen},\R}(A) := \frac{\sqrt{2}\|A\|_F}{\dist(A,\Sigma)}.
\end{equation}
This definition is motivated by a result by Wilkinson~\cite{Wilkinson72}, 
which (in slight reformulation) states that for all real eigenvalues $\lambda$ of $A$
\begin{equation}\label{eq:wilki}
 \kappa(A,\lambda) \le \frac{\sqrt{2}\|A\|_F}{\dist(A,\Sigma)} .
\end{equation}
Thus $\kappa_{{\sf eigen},\R}$ is a conic condition number;  
it varies continuously with $A$ and takes the value $\infty$ when 
$A\in\Sigma$.

In~\cite{Demmel88}, Demmel used the fact that 
the right-hand side of (\ref{eq:wilki}) is conic 
to obtain bounds on the
tail of $\kappa_{{\sf eigen},\R}(A)$ for random $A$. We next
use this same fact to obtain smoothed analysis estimates.

\begin{proposition}\label{prop:eigen}
For all $n\geq 1$ and $0<\sigma\le 1$ we have
$$
  \sup_{A\in S^{n^2-1}}
   \bfE_{Z\in B(A,\sigma)}(\ln \kappa_{{\sf eigen},\R}(Z))\leq
    8\ln n+2\ln\frac1\sigma + 6.
$$
\end{proposition}

\begin{Proof}
Let $W$ be the set of matrices having multiple eigenvalues 
(real or complex). Clearly, $\Sigma\subseteq W$. In addition, 
$W$ is the zero set of the discriminant polynomial, 
which is a homogeneous polynomial of degree $n^2-n$ 
in the entries of~$A$ (compare~\cite[Prop.~3.4]{BCL:06a}).
Theorem~\ref{th:main} applied to the conic condition number
$\frac{\|A\|_F}{\dist(A,\Sigma)}$ implies the stated bound.
\end{Proof}

\subsection{General Eigenvalue computations}\label{ss:g_eigen}

In addition to the computation of real eigenvalues, one can 
consider the computation of general (possibly complex) eigenvalues 
for real (or complex) matrices $A$. Replacing the inner product 
in $\R^n$ by a corresponding Hermitian product in $\C^n$, the 
derivation of~(\ref{eq:cnlambda}) still holds and one can now 
define  
$$
   \kappa_{\sf eigen}(A) := \max_{\lambda}\kappa(A,\lambda)
$$
where  the maximum is over all eigenvalues $\lambda$ of $A$. 
Replacing now $\Sigma$ by the set $W$ in the proof of 
Proposition~\ref{prop:eigen}, Wilkinson's result~\cite{Wilkinson72} 
still holds and yields 
$$
  \kappa_{\sf eigen}(A)\leq \frac{\sqrt{2}\|A\|_F}{\dist(A,W)}
$$
from which the following result easily follows (for part (ii) 
one considers complex numbers as pairs of real numbers, i.e., 
a complex matrix $A$ as $A\in\R^{2n^2}$; this is actually 
the way programming languages deal with complex numbers). 

\begin{proposition}\label{prop:g_eigen}
For all $n\geq 1$ and $0<\sigma\le 1$ we have
\begin{description}
\item[(i)]
For real matrices $A$,
$$
  \sup_{A\in S^{n^2-1}}
   \bfE_{Z\in B(A,\sigma)}(\ln \kappa_{\sf eigen}(Z))\leq
    8\ln n+2\ln\frac1\sigma + 6.
$$
\item[(ii)]
For complex matrices $A$,
$$
  \sup_{A\in S^{2n^2-1}}
   \bfE_{Z\in B(A,\sigma)}(\ln \kappa_{\sf eigen}(Z))\leq
    8\ln n+2\ln\frac1\sigma + (6+2\ln 2).
$$
\end{description}
\end{proposition}

\begin{Proof}
Part (i) is as in Proposition~\ref{prop:g_eigen}. For part~(ii) 
one notes that $W$, as a subset of $\R^{2n^2}$, is the zero set of 
the real and imaginary parts of the discriminant polynomial, which have 
both degree $n^2-n$. 
\end{Proof}

%%%%%%%%%%%%%%%%%%%%%%%%%%%%%%%%%%%%%%%%%%%%%%%%%%%%%%%%
\subsection{Solving real polynomial systems}\label{ss:systems}

Let $d_1,\ldots,d_n\in\N\setminus\{0\}$. We denote by $\Hd$ the
vector space of polynomial systems $f=(f_1,\ldots,f_n)$ with
$f_i\in\R[X_0,\ldots,X_n]$ homogeneous of degree $d_i$.
For $f,g\in \Hd$ we write
\begin{equation*}
   f_i(x) = \sum_{\alpha} a^i_\alpha X^\alpha, \quad g_i(x)
   = \sum_{\alpha} b^i_\alpha X^\alpha,
\end{equation*}
where $\alpha=(\alpha_0, \dots, \alpha_n)$ is assumed to range over
all multi-indices such that $|\alpha| = \sum_{k=0}^n \alpha_k = d_i$
and $X^\alpha:= X_0^{\alpha_0}X_1^{\alpha_1}\cdots X_n^{\alpha_n}$.
We endow the space $\Hd$ with the inner product
$\langle f,g \rangle := \sum_{i=1}^n \langle f_i, g_i \rangle$, where
\begin{equation*}
  \langle f_i, g_i \rangle =
  \sum_{|\alpha|=d_i}
   {d_i \choose \alpha}^{-1}\, a^i_\alpha\, b^i_\alpha.
\end{equation*}
Hereby, the multinomial coefficients are defined as
\begin{equation*}
  {d_i\choose \alpha} = \frac{d_i!}
  {\alpha_0 ! \alpha_1 ! \cdots \alpha_n !}.
\end{equation*}
This inner product has the beautiful property of being invariant
under the natural action of the orthogonal group $O(n+1)$ on $\Hd$.
In the case of one variable, it was introduced by Weyl~\cite{weyl:50}.
Its use in computational mathematics goes back at
least to Kostlan~\cite{Kostlan93}.

Throughout this section, let $\|f\|$ denote the corresponding norm 
of $f\in\Hd$.
The Weyl inner product defines a Riemannian structure on the sphere
$S(\Hd):=\{f\in\Hd \mid \|f\|=1\}$. As in the previous sections,
we endow this sphere with the corresponding projective distance $d_{\P}$.

In a seminal series of papers, M.~Shub and
S.~Smale~\cite{Bez1,Bez2,Bez3,Bez4,Bez5} studied the problem of,
given $f\in\Hd\otimes_\R\C$, to compute an approximation of a 
complex zero of~$f$.
They proposed an algorithm and studied its complexity in terms of, among
other parameters, a condition number $\munorm(f)$ for $f$.
In the following we will recall the definition of 
$\munorm(f)$ adapted to the case of real systems and real zeros  
(see~\cite[Chapter~12]{bcss:95} for details). 

For a simple zero $\z\in S^n$ of $f\in\Hd$ one defines
\begin{equation*}
  \munorm(f,\z):= \|f\|\left\| (Df(\z)_{|T_\z})^{-1}
  \mathsf{diag}(\sqrt{d_1},\ldots,\sqrt{d_n})\right\|,
\end{equation*}
where $Df(\z)_{|T_\z}$ denotes the restriction of the derivative of
$f\colon\R^{n+1}\to\R^n$ at $\z$ to the tangent space
$T_\z S^n =\{v\in\R^{n+1}\mid \langle v,\z\rangle = 0\}$
of $S^n$ at $\z$.
(The norm on the right is the operator norm with respect to the 
Euclidean norm.)
If $\z$ is not a simple root of $f$ we set $\munorm(f)=\infty$.
Note that $\munorm(f,\z)$ is homogeneous of degree~0 in $f$ and $\z$.

Shub and Smale~\cite{Bez2} proved a condition number theorem
for the condition number $\munorm(f,\z)$ for complex polynomial
systems $f$ and complex roots $\z$.
To describe a corresponding result in the real situation, 
consider for $\z\in S^n$
$$
 \Sigma_\z:= \{g\in \Hd \mid \mbox{$\z$ is a multiple zero of $g$}\} .
$$

\begin{theorem}\label{th:CNTR}
For a zero $\z\in S^n$ of $f\in S(\Hd)$ we have
\begin{equation*}
   \munorm(f,\z)=\frac{1}{d_{\P}(f,\Sigma_{\z}\cap S(\Hd))}.
\end{equation*}
\end{theorem}

\begin{proof}
The proof of~\cite{Bez2} (see also \cite[\S12.4]{bcss:95}) carries 
over immediately to the real situation.
\end{proof}

Let $\Sigma\subseteq\Hd$ be the set of systems 
$f\in \Hd$ having a real multiple zero. 
Note that 
$\Sigma_\z\subseteq\Sigma$. Therefore, 
by Theorem~\ref{th:CNTR}, 
\begin{equation*}\label{eq:ss}
  \munorm(f,\z) \le \frac1{d_\P(f,\Sigma\cap S(\Hd))} .
\end{equation*}
We define the condition number for real polynomial solving as 
the right-hand side: 
$$
   \munormR(f) := \frac1{d_\P(f,\Sigma\cap S(\Hd))} .
$$
This is by definition a conic condition number: 
it varies continuously with $f$ and takes the value $\infty$ 
when $f\in\Sigma$.

\begin{proposition}\label{cor:realpolysolv}
For all $d_1,\ldots,d_n\in\N\setminus\{0\}$ and all $\sigma\in(0,1]$
we have
\begin{equation*}
 \sup_{f\in S(\Hd)}\ \bfE_{g\in B(f,\sigma)} (\ln\munormR(g))\leq
   2\ln N + 4\ln \DD + 2\ln n + 2\ln\frac1{\s} + 7,
\end{equation*}
where $N=\dim\Hd -1$ and $\DD=d_1\cdots d_n$ is the B\'ezout number.
\end{proposition}

\begin{Proof}
Let $W$ be the {\em discriminant variety} consisting of the systems 
$f\in \Hd$ having a complex multiple zero. Then $\Sigma\subseteq W$. 
In addition, it is well known that $W$ is the 
zero set of a multihomogeneous polynomial of
total degree bounded by $2n\DD^2$, where
$\DD= d_1\cdots d_n$ is the B\'ezout number
(see, e.g., \cite[Lemma 3.6]{BCL:06a}).
The statement now follows immediately from 
Theorem~\ref{th:main}.
\end{Proof}

\begin{remark}
\begin{description}
\item[(i)]
By the results in~\cite{Bez1}, the condition number 
$\munorm(f)$ not only measures the maximum sensitivity of complex 
solutions to the input~$f$, but it is also a crucial complexity 
parameter for algorithms approximating such solutions. 
While $\munormR(f)$ shares with $\munorm(f)$ the first property, 
it is not clear whether this is also the case for the second one.
\item[(ii)]
By considering complex numbers as pairs of real numbers one can 
see the problem of, given $f\in\Hd$ (or in $\Hd\otimes\C$), compute 
an approximation of a complex zero of $f$ (studied by 
Shub and Smale~\cite{Bez1,Bez2,Bez3,Bez4,Bez5}) as a problem 
over the reals. Proceeding as in \S\ref{ss:g_eigen}, one can find 
bounds for the smoothed analysis of this problem similar 
to those in Proposition~\ref{cor:realpolysolv}. 
\end{description}
\end{remark}

%%%%%%%%%%%%%%%%%%%%%%%%%%%%%%%%%%%%%
\subsection{Real zero counting}

Consider the problem of, given $f\in\Hd$, counting the number of 
real zeros of $f$ in~$S^n$. Unlike the problems considered so far, 
this is a problem with a discrete output. This means that 
sensitivity considerations as described in the previous problems 
do not apply here. Yet, finite precision algorithms will require 
more precision to give a reliable output 
when the input $f$ is close to the set 
$\Sigma$ of systems with multiple real zeros, and will 
not give any such reliable output when $f\in\Sigma$ (since 
in this case, arbitrarilly small perturbations of $f$ will 
change the output). It therefore makes sense to define 
as condition number for the counting problem
$$
   \cond(f)=\frac{\|f\|}{\dist(f,\Sigma)} .
$$
Since $\cond(f) = \munormR(f)$ for $f\in S(\Hd)$, 
the bounds of Proposition~\ref{cor:realpolysolv} hold for 
$\cond(f)$ as well. We note, however, that in 
\S\ref{ss:systems} we assumed a specific inner product 
on $\Hd$ (whose properties are crucial in the proof 
of Theorem~\ref{th:CNTR}). In contrast, for the smoothed analysis 
of $\cond(f)$, any inner product on $\Hd$ would do. 

%%%%%%%%%%%%%%%%%%%%%%%%%%%%%%%%%%%%%%%%%%%%%%%%%%%%
\section*{Appendix}

\setcounter{section}{1} \setcounter{subsection}{0}
\setcounter{theorem}{0}
\renewcommand{\thesection}{\Alph{section}}

The constants $\kc(p,i)$ appearing in Theorem~\ref{th:kin-form-S}
depend only on $p$ and $i$ and are independent of the manifold $M$.
We derive the expression for $\kc(p,i)$ stated in 
Theorem~\ref{th:kin-form-S} by selecting a simple enough $M$.
Consider the boundary $M_\a=M_\a(q)$ of the ball $B_R(q,\a)$ in
$S^p$ of radius $0<\a\le\pi/2$ centered at~$q$ (recall
Example~\ref{ex:geoball}). According to Theorem~\ref{th:kin-form-S}
we have for $0\le i <p-1$
\begin{equation*}
 \mu_i(M_\a) = \kc(p,i) \frac1{\Oh_p}\int_{S^p}
 \mu_i(M_\a(z)\cap S^{i+1}) dS^p(z) .
\end{equation*}
Let $\rho(z)$ denote the Riemannian distance from $z\in S^p$ to
$S^{i+1}$. If $\rho(z)<\a$, then the intersection $M_\a(z)\cap
S^{i+1}$ is the boundary of a ball in $S^{i+1}$. The radius $\d(z)$ of
the sphere $M_\a(z)\cap S^{i+1}$ satisfies $\cos\a = \cos\rho(z)\cdot
\cos\d(z)$ by a well known formula of spherical trigonometry.  From
Example~\ref{ex:geoball} we know that $\mu_i(M_\a(z)\cap S^{i+1}) =
\Oh_i\, (\cos\d(z))^i$.  On the other hand, $M_\a(z)$ does not
intersect $S^{i+1}$ if $\rho(z) >\a$.

>From these reasonings we obtain
\begin{equation*}
 \frac1{\Oh_p}\int_{S^p} \mu_i(M_\a(z)\cap S^{i+1}) dS^p(z)
 = \frac{\Oh_i}{\Oh_p}\int_{0}^\a \left(\frac{\cos\a}{\cos\rho}\right)^i
 \frac{d}{d\rho}\,\vol_p T_{\P}(S^{i+1},\rho)\, d\rho .
\end{equation*}
>From Lemma~\ref{le:J1} and the definition of $J_{p,p-i-1}(\rho)$
it follows that
\begin{equation*}
 \frac{d}{d\rho}\,\vol_p T_R(S^{i+1},\rho) =
 \Oh_{i+1} \Oh_{p-i-2} (\sin\rho)^{p-i-2}\,(\cos\rho)^{i+1}.
\end{equation*}
We obtain
\begin{eqnarray*}
 \frac1{\Oh_p}\int_{S^p} \mu_i(M_\a(z)\cap S^{i+1}) dS^p(z) &=&
 \frac{\Oh_i\Oh_{i+1} \Oh_{p-i-2}}{\Oh_p}\int_{0}^\a (\cos\a)^i
    (\sin\rho)^{p-i-2}\,\cos\rho\,d\rho.\\
 &=& \frac{\Oh_i\Oh_{i+1} \Oh_{p-i-2}}{\Oh_p} (\cos\a)^i\
    \frac{(\sin\a)^{p-i-1}}{p-i-1}.
\end{eqnarray*}
On the other hand, by Equation~(\ref{eq:muiM}),
\begin{equation*}
 \mu_i(M_\a) = {p-1\choose i} \Oh_{p-1} (\sin\a)^{p-i-1} (\cos\a)^{i},
\end{equation*}
By comparing the last two equations, the asserted form of $\kc(p,i)$
follows.
\hfill$\Box$

{\small
%\bibliography{condition}
%\bibliography{../biblio/lit-bank}
%\input hilbert-rev.bbl

}
\end{document}